# Optimizing Traffic Signal Control for Continuous-Flow Intersections: Benchmarking against a State-of-Practice Model

Yining Hu[1], David Rey[2], Reza Mohajerpoor[3], Meead Saberi[1,*]

[1] *Research Centre for Integrated Transport Innovation (rCITI), School of Civil and Environmental Engineering, University of New South Wales (UNSW), Sydney, Australia, 1466*
[2] *SKEMA Business School, Université Côte d'Azur, Sophia Antipolis, France*
[3] *School of Civil Engineering, The University of Sydney, Sydney, Australia*

[*] *Corresponding author (meead.saberi@unsw.edu.au)*

**Abstract**: Continuous-Flow Intersections (CFI), also known as Displaced Left-Turn (DLT) intersections, aim to improve the efficiency and safety of traffic junctions. A CFI introduces additional sub-intersections upstream of the main intersection to split the left-turn flow from the through movement before it arrives at the main intersection, decreasing the number of conflict points between left-turn and through movements. This study develops and examines a two-step optimization model for CFI traffic signal control design and demonstrates its performance across more than 300 different travel demand scenarios. The proposed model is compared against a benchmark state-of-practice CFI signal control model. Microsimulation results suggest that the proposed model reduces average delay by 17% and average queue length by 32% for a full CFI compared with the benchmark signal control model.

## INTRODUCTION

Continuous-Flow Intersections (CFI), also known as Displaced Left-Turn (DLT) intersections, have been proposed to enhance the safety and efficiency of intersections with a high volume of left-turn traffic (Al-Salman and Salter, 1974; Hutchinson, 1995; Goldblatt et al., 1994; Zhao et al., 2018). A CFI employs distinct sub-intersections to divert the left-turn flow from the through movement before it reaches the main intersection, reducing the number of phases at a signalized intersection from four to two. Left-turning vehicles pass through a sub-intersection first, moving to the left side of the road before arriving at the main intersection, thereby eliminating conflict points between left-turn and through movements (Simmonite and Chick, 2004; Steyn et al., 2014). This design allows left-turn and through traffic to share the same phase, enhancing traffic efficiency and safety. CFIs are particularly suitable for intersections with heavy left-turn traffic (Reid and Hummer, 2001).

While the traffic signal control problem of CFIs has been investigated in past studies (Suh and Hunter, 2014; Shokry et al., 2018; Zhao et al., 2015), existing literature has focused on a limited number of scenarios without benchmarking against other traffic signal control design methods. Traffic demand at an intersection can follow complex patterns, impacting the actual performance of CFI implementation. In this study, we adopt and modify a two-step optimization model for CFI traffic signal control design, originally proposed by Yang and Cheng (2017). We simulate and test the performance of the proposed model on over 300 different demand scenarios, including various Origin-Destination (OD) demand patterns and different demand levels. Comparisons are made between the developed model and a state-of-practice CFI traffic signal control design model. Extensive microsimulation results suggest that, on average, the proposed model exhibits improved performance compared to the benchmark model, particularly for intersections with left-right-turn traffic.

This study makes several key contributions. Firstly, it adapts and enhances the traffic signal optimization model proposed by Yang and Chen (2017) specifically for application to a full CFI. Secondly, it conducts a rigorous comparison between the developed model and a state-of-practice benchmark model derived from a Federal Highway Administration (FHWA) study. Through extensive simulation-based evaluation encompassing over diverse demand scenarios, the study thoroughly examines the performance of the proposed optimization model against the benchmark. The results demonstrate consistently superior performance of the developed optimization model across a wide range of demand scenarios, highlighting its effectiveness in enhancing traffic efficiency and indicating its potential for real-world application.

Note that this study explores the CFI traffic signal control problem in Australia, a left-hand side driving country. Therefore, in the methodology and numerical results sections of the manuscript, the modeling framework and simulation outcomes are based on left-hand side driving rules in which the CFI includes a displaced right turn as opposed to displaced left turn.

**LITERATURE REVIEW**
The concept of CFI was first proposed a few decades ago (Al-Salman and Salter, 1974) and has been explored in numerous studies since then. The literature already includes several previous research studies aiming to optimize the traffic signal control design of CFIs using different approaches. Goldblatt et al. (1994) evaluated the performance of a CFI using the TRAF-NETSIM microscopic simulation model against a conventional intersection. They showed that CFIs offer many advantages over conventional intersections in terms of traffic performance. In a later study, Reid and Hummer (2001) showed that CFIs outperform conventional intersections, particularly when heavy left-turn volumes exist. Several other research studies also revealed that CFIs can indeed improve traffic efficiency compared to conventional intersections (Jagannathan and Bared, 2004; Hildebrand, 2007; Esawey and Sayed, 2007; Dhatrak et al., 2010; Park and Rakha, 2010; Olarte et al., 2011; Autey et al., 2013; Zhao et al., 2015; Pan et al., 2021).

*Safety*
The safety aspect of CFIs has also been investigated by several studies. Theoretically, a CFI contains a fewer number of conflict points compared to a conventional intersection, thereby reducing rear-end crash rates (Steyn et al., 2014). However, past research has also revealed potential risks associated with CFIs. Yahl (2013) conducted an observational before-and-after study to investigate the safety impact of CFIs at five sites in the U.S. The results showed that only one out of the five studied sites experienced a decrease in the number of collisions, while the other sites exhibited slightly decreased safety performance due to variations in their designs. Similar conclusions were also reached by Abdelrahman et al. (2020). Qu et al. (2020) conducted a case study over two CFIs in Texas. The results indicated that the implementation of CFIs did not increase the overall frequency of crashes at the studied intersections and reduced left-turn-related collisions at the main intersection, as well as rear-end crashes at the right-turn merge points. Many crashes occurred due to illegal left-turns and drivers not being familiar with the new geometric design. Zlatkovic (2015) also investigated the safety performance of CFIs and demonstrated that they can improve traffic safety.

*Signal Control Optimization*
Suh and Hunter (2014) developed a two-model framework for the CFI signal control problem. The first model optimizes the green bandwidth, representing the time gap that vehicles can pass two consecutive intersections without being stopped by the red light (Shokry et al., 2018). The second model minimizes intersection delay to improve traffic efficiency. Zhao et al. (2015) proposed a model that optimizes the intersection design in addition to CFI signal control. They used two types of evaluation measures. The first measure is called the common flow multiplier, indicating the ratio between actual capacity and traffic volume, while the second measure represents the bandwidth. Yang and Cheng (2017) developed two complementary optimization models for an asymmetric two-leg CFI signal control problem. The first model follows a two-step sequential approach that uses the common flow multiplier in the first step and the bandwidth in the second step. Meanwhile, the second model merged the two steps into a single multi-objective problem. In a later study, Shokry et al. (2018) studied the optimization of the offset between a series of CFIs to maximize the bandwidth. A few studies considered pedestrians and cyclists. Zhao et al. (2019) proposed an optimisation design approach for left-turning bicycles to improve the practical capacity of CFI. Jiang and Gao (2020) considered constraints related to safe-crossing requirements of pedestrians and non-motorized vehicles in the CFI traffic signal control problem.

Common flow multiplier and bandwidth are two widely used performance measures in traffic signal optimization models in the literature especially those based on optimizing the progression for multi-intersection coordination ( Shirvani and Maleki, 2016; Lu et al, 2017; Shorky et al., 2018; Xu et al., 2019; Ma et al., 2020). Common flow multiplier is directly related to intersection capacity, which can be linked with intersection efficiency. While bandwidth represents efficiency and progression in a CFI, it may not necessarily provide a direct link with intersection delay (Suh and Hunter, 2014). Bandwidth establishes a connection between the main intersection and sub-intersections in a CFI, enabling the optimization of offsets.

In a conventional intersection signal control optimization problem, three main variables are often considered, including cycle length, phase design, and signal timing. For a standard symmetric CFI, the signal phasing problem is less relevant and significant. Instead, the offset is considered one of the key variables. The model developed by Suh and Hunter (2014) optimizes the through movement phase time for secondary intersections and the offset for each secondary intersection. Zhao et al. (2015) proposed a model that includes various decision variables. Besides the signal timing for the main intersection and sub-intersections, the model also optimizes the selection of intersection type, the length of the displaced left-turn lane, and the lane markings. In a typical CFI, left-turn movements are separated from the main intersection, and each of the five intersections in a CFI includes only two phases. In another study, Yang and Cheng (2017) investigated the signal control optimization problem for an asymmetric two-leg CFI. The presence of conventional intersection legs necessitates multiple phases at the main intersection, reintroducing the phasing problem eliminated in the full CFI. To incorporate signal phasing into the optimization, the authors separated the eight movements and identified five critical paths instead of predetermined phases.

In a CFI, the sub-intersections are usually situated close to the main intersection, requiring vehicles to traverse two intersections within a short distance. Given the limited capacity of turning bays between the sub-intersections and the main intersection, queue spillbacks are common if progressions between intersections are not well-coordinated or planned. Consequently, the offset between the main intersection and sub-intersection becomes a critical variable in the CFI traffic signal design problem. While existing literature has investigated various objective functions and decision variables in the context of CFI traffic signal control, most studies have only tested the performance of proposed signal control plans with a single demand scenario as the base case (Yang and Cheng, 2017; Shokry et al., 2018), or with a limited number of demand scenarios (Suh and Hunter, 2014; Zhao et al., 2015). In this study, we investigate the CFI traffic signal control problem in Australia, a left-hand driving country. Accordingly, the term "left-turn" in most other studies is equivalent to "right-turn" in this paper. The study proposes a two-step traffic signal control optimization model for a complete symmetric CFI and evaluates the performance of the model under several demand scenarios, covering a variety of traffic patterns across all movements. Existing studies have rarely compared the developed models with other CFI signal optimization methods, especially those used in practice. Most studies have focused on comparing CFIs with conventional intersections (Zhao et al., 2015; Shokry et al., 2018), or comparing the developed models with simulation-based signal plans (Yang and Cheng, 2017; Qu et al., 2021). In this study, we provide a comprehensive comparative analysis between the proposed model and the state-of-the-practice model proposed by Qi et al. (2020) and Qu et al. (2021), referred to as the benchmark model. Microsimulation results demonstrate that the proposed model outperforms the benchmark model in most scenarios, particularly in scenarios with heavy right-turn traffic.

**MATHEMATICAL FORMULATION**
We formulate a mathematical model based on the two-step model originally proposed by Yang and Cheng (2017). The original model was developed for an asymmetric CFI with two CFI legs and two conventional legs. Consequently, more than two phases are needed at the main intersection for the right-turn movements from the conventional legs, thus requiring Yang and Cheng (2017) to consider the phasing problem as well. While many other CFI signal design models adopt a pre-determined two-phase model in optimization, the model proposed by Yang and Cheng (2017) offers more flexibility to extend to more sophisticated cases. For example, pedestrians and cyclists can be included in the model by introducing additional phases at the intersection (Li et al., 2010; Ma et al., 2015). We extend the original two-leg asymmetric CFI proposed by Yang and Cheng (2017) to a fully symmetric CFI. Constraints and variables related to the conventional legs are adjusted to accommodate the full CFI (see Fig. 2). Constraints for the green band in the base model did not cover all possible coordination situations between the intersections. Therefore, some available green time was not fully utilized in the green band. To address this limitation, we formulated new constraints so that longer green times can be covered when optimizing the bandwidth. The developed mathematical model consists of two steps. The first step determines a common cycle length for the CFI, as well as the green time for each movement at each intersection. The second step determines the offset at each intersection.

For a complete list of notations, please refer to the Notation section at the end of the manuscript.

*Step 1*

Like other coordinated intersections in a traffic network (H. Qi et al., 2020), the main intersection and sub-intersections in a CFI share a common cycle length. In the first step of the optimization, the common cycle length for the CFI and the green split for each movement at each intersection are determined. To establish a linear optimization model, the reciprocal of cycle length $\xi$ and the reciprocal of green time $\phi_{l,i}$ are used as decision variables. The objective function is to maximise the sum of multipliers $\mu_l$ that represent the capacity at the intersection, as expressed in Equation (1). Equation (2) ensures the degree of saturation with the multiplier $\mu_l$ does not exceed the physical capacity.

$$\text{Maximise} \sum_{l \in L} \mu_l \tag{1}$$

$$\mu_l \alpha_{l,i} Q_i \leq s_{l,i}(\phi_{l,i} - \delta \times \xi) \quad \forall l \in L_i, i \in I \tag{2}$$

The additional constraints for maximum and minimum cycle length and green time are as expressed below:

$$\frac{1}{C_{\max}} \leq \xi \leq \frac{1}{C_{\min}} \tag{3}$$

$$\xi \times g_{\min,i} \leq \phi_{l,i} \leq \min(\xi \times g_{\max,i}, 1) \quad \forall l \in L_i, i \in I \tag{4}$$

To ensure movements are assigned to the phase with no conflict at each intersection, several constraints are defined that are summarized in Appendix A.

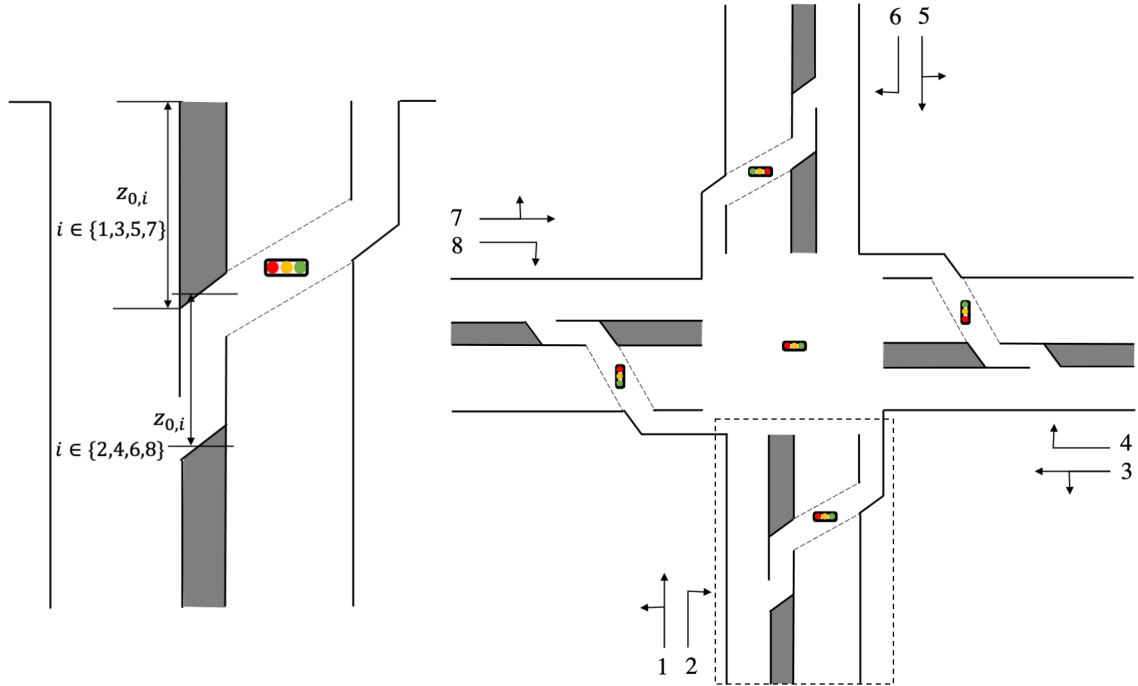

**Fig. 1.** Maximum acceptable distance between stop line and end of queue for through movements at the main intersection and right-turn movements at sub-intersections in a CFI.

To avoid queue spillback at the right-turn bays and the blockage of the right-turn bays caused by the through movement queues, the maximum distance between stop line and rear of queue $\tau_i$ of right-turn flows at sub-intersections and the queue of through movements at the main intersection, as estimated in Equation (5), should not exceed the corresponding capacity $z_{0,i}$ (See Equation (6)). The multiplicand in Equation (5) calculates the maximum number of vehicles in the queue, and together with the multiplier, the expression calculates the maximum distance from the stop line to the end of the queue. $z_{0,i}$ is the length of the right-turn bay for right-turn movements that prevents the spillback from the bay to occur as illustrated in Fig. 1, while for through movements, $z_{0,i}$ is the distance between the main intersection and the entrance of the right-turn bay to prevent the through movements to block the right-turn Bay.

$$\tau_i = \frac{s_{f,i}}{s_{f,i} - \alpha_{f,i}Q_i} \times \frac{(1 - \phi_{f,i} + \delta \times \xi) \cdot \alpha_{f,i}Q_i}{\xi} \qquad \forall i \in I \qquad (5)$$

$$\tau_i \leq z_{0,i} \qquad \forall i \in I \qquad (6)$$

where $f$ is the first intersection passed by movement $i$.

Equation (7) is formulated to prevent the intersection from oversaturation. On the left side of the equation, the nominator calculates the maximum number of vehicles in the queue, the denominator calculates the dissipation speed of the queue, and the fraction indicates the duration of the queue dissipation, which should not be larger than green time. Since we are assuming a uniform distribution for the arriving traffic, the average number of vehicles that arrive within a cycle should not exceed the maximum number of vehicles that can be served within the green time, otherwise the queue accumulates and eventually spills back.

$$\frac{(1 - \phi_{f,i} + \delta \times \xi) \cdot \alpha_{f,i}Q_i}{(s_{f,i} - \alpha_{f,i}Q_i)\xi} \leq \frac{\phi_{f,i}}{\xi} \qquad \forall i \in I \qquad (7)$$

After obtaining the optimal values of $\xi$ and $\phi_{l,i}$ (inclusive of $\phi_{f,i}$ for the first intersection passed by the movements) the cycle length and the green splits can be determined by taking the reciprocal of $\xi$ and $\phi_{l,i}$ as expressed in Equation (8) and (9), which will be used as input parameters in the second step of the optimization.

$$C = \frac{1}{\xi} \qquad (8)$$

$$g_{l,i} = \phi_{l,i}C \qquad \forall l \in L_i, i \in I \qquad (9)$$

***Step 2***
In this step, we determine the offsets of each intersection's signal $\theta_l$. The objective function in this step aims to maximize the sum of the weighted green bandwidth $b_i$ that every movement $i$ can receive:

$$\text{Maximise} \sum_{i \in I} \eta_i b_i \qquad (10)$$

The following constraints are defined to ensure the solution feasibility.

$$w_{l,i} \geq 0 \qquad (11)$$

$$w_{l,i} + b_i \leq g_{l,i} \qquad \forall l \in L_i, i \in I \qquad (12)$$

where $w_{l,i}$ is the time between the start of a green phase $g_{l,i}$ and the start of the corresponding green band $b_i$.

We have identified six different scenarios for the green band between two consecutive intersections (See Fig. 2), in which cases 1 to 5 are scenarios that Movement $a$ receives an effective green band, while in case 6 it receives no green band. In case 1, the green band starts later than the green phase and ends together with it. In case 2, the green band starts together with the green phase and ends earlier than it. In case 3, the green band is separated into two parts by the red phase of the second intersection. In case 4, the green band starts later and ends earlier than the green phase. In case 5, the green band covers the whole green phase. In case 6, no overlapping period exists between green phases of the two intersections. Corresponding constraints are provided in Appendix C.

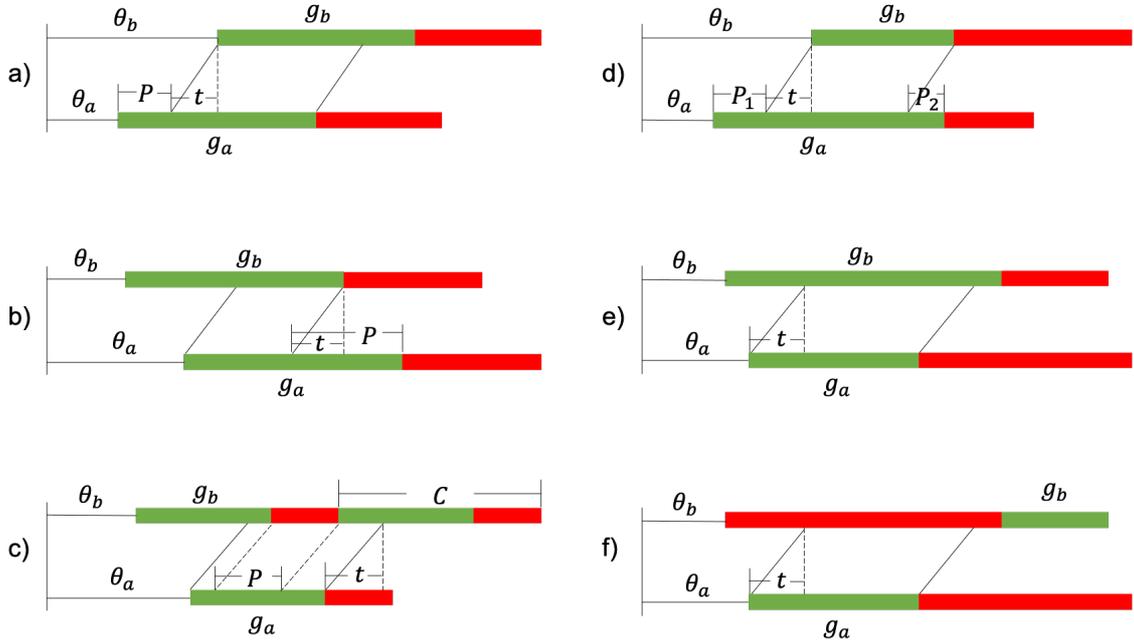

**Fig. 2.** Different possible green band progression scenarios between intersection $a$ and intersection $b$: a) case 1, b) case 2, c) case 3, d) case 4, e) case 5, and f) case 6.

The sub-intersections in a CFI often have a short distance to the main intersection. To avoid spillback from the main intersection to the sub-intersections, we introduce a new variable, uncoordinated time period $P_{d,i,k}$, which will be explained in detail later in this section. The green bandwidth is the time period that a vehicle can travel through multiple intersections without a stop, which should not be larger than the difference between the green time and the uncoordinated time period $P_{d,i,k}$ between a pair of intersections $d = (a, b)$. $k$ is an index for the 6 green band cases in Fig. 2, $P_{d,i,k}$ is 0 if case $k$ does not apply for intersection $i$ between intersection pair $d$.

$$b_i \leq g_{a,i} - P_{d,i,k} \quad \forall d \in D_i, (a,b) \in d, i \in I, k \in K \quad (13)$$

A binary variable $y_{d,i,k}$ is introduced to allow the model to select the proper type of green-band from Cases 1 to 6 (as shown in Fig. 2) for each movement $i$:

$$y_{d,i,k} = \begin{cases} 0; & \text{if movement } i \text{ obtains case } k \text{ green band between intersection pair } d \\ 1; & o.w. \end{cases} \quad (14)$$

For movements with green band progression, the bandwidth should not be smaller than the minimum effective bandwidth $b_e$ as the constraint given below, $M$ is a very large positive value.

$$b_i \geq b_e - (1 - y_{d,i,6})M \quad \forall d \in D_i, i \in I \quad (15)$$

The following constraints enforce the progression between intersections:

$$\theta_a + r_{a,i} + w_{a,i} + t_{a,i} + n_{a,i}C \geq \theta_b + r_{b,i} + w_{b,i} + n_{b,i}C - (1 - y_{d,i,6})M \\ \forall d \in D_i, i \in I \quad (16)$$

$$\theta_a + r_{a,i} + w_{a,i} + t_{a,i} + n_{a,i}C \leq \theta_b + r_{b,i} + w_{b,i} + n_{b,i}C + (1 - y_{d,i,6})M \\ \forall d \in D_i, i \in I \quad (17)$$

Let $d \in D_i$ be a pair of intersections and let $(a, b) \in d$ be the intersections of pair $d$, representing the first and second intersections passed by movement $i$ in the intersection pair $d$. The uncoordinated time period $P_{d,i,k}$ (as also shown in Fig. 2) of movement $i$ between intersections $a$ and $b$ (which are in intersection pair $d$) in scenario $k$ can be calculated using Equations (18) to (23):

$$P_{d,i,1} = g_{a,i} - (\theta_a + r_{a,i} + n_{a,i} \times C + g_{a,i} - \theta_b - r_{b,i} - n_{b,i} \times C + t_{d,i}) \\ \forall d \in D_i, (a,b) \in d, i \in I \quad (18)$$

$$P_{d,i,2} = g_{a,i} - (\theta_b + r_{b,i} + n_{b,i} \times C + g_{b,i} - \theta_a - r_{a,i} - n_{a,i} \times C - t_{d,i}) \\ \forall d \in D_i, (a,b) \in d, i \in I \quad (19)$$

$$P_{d,i,3} = C - g_{b,i} \quad \forall d \in D_i, (a,b) \in d, i \in I \quad (20)$$

$$P_{d,i,4} = g_{a,i} - g_{b,i} \quad \forall d \in D_i, (a,b) \in d, i \in I \quad (21)$$

$$P_{d,i,5} = 0 \quad \forall d \in D_i, (a,b) \in d, i \in I \quad (22)$$

$$P_{d,i,6} = g_{a,i} \quad \forall d \in D_i, (a,b) \in d, i \in I \quad (23)$$

To prevent spillback between the main and sub intersections, the below constraint ensures the queue accumulated within a cycle should not be longer than the capacity of the bay between the two intersections:

$$Q_i \times C \times \frac{P_{d,i,k}}{g_{a,i}} \leq z_{d,i} + y_{d,i,k} \times M \quad \forall d \in D_i, i \in I, k \in K \quad (24)$$

Similar to $\phi_{l,i}$, constraints for $r_{l,i}$ are defined in Appendix B.

The proposed model is a mixed integer linear program (MILP) and can be solved by commercial optimization solvers.

**MODEL EVALUATION AND BENCHMARKING**
We evaluate the performance of the proposed optimization model for a full symmetric CFI with sub-intersections on all four legs. The model is programmed in Python, and the GNU Linear Programming Kit (GLPK) is used as the optimiser, with the solver time limit set to 20 seconds. The intersection is tested with more than 300 different traffic demand scenarios, and the corresponding microsimulations are conducted in AIMSUN with the signal plan generated by the proposed optimization model. The simulation duration is set to 1 hour with a fixed signal control. Each simulation contains 5 replications. The average delay and queue length in the experiment results are directly obtained from AIMSUN.

*Benchmark model*
A state-of-practice signal timing strategy and model for CFIs is proposed by Qi et al. (2020) as part of a U.S. Federal Highway Administration (FHWA) study, in which a case study was conducted. The results showed that the developed model outperformed Synchro default signal timing and Synchro optimized signal timing. We select this model as a state-of-practice benchmark and evaluate the performance of our proposed model against it. Unlike our proposed model, which provides a different signal plan for each sub-intersection, the benchmark model in Qi et al. (2020) generates a symmetric signal design with identical phasing and timing for two sub-intersections on the opposite side of the main intersection. In the benchmark model, the signal timing of the main intersection is determined by the ratios of traffic volumes. The timing of sub-intersections' signals as well as the offsets are determined based on the traffic progression between the intersections.

*Scenario design*
281 out of 322 designed demand scenarios include a fixed total travel demand of 4,000 vehicles per hour. According to preliminary experiments, 4,000 vehicles per hour represent a level of demand that is neither too high, leading to oversaturation in some sections, nor too low, indicating that intersection efficiency improvement is unnecessary. In each scenario, the total demand is allocated differently among the 8 movements, ranging from 3% to 34% of the total demand. The remaining 41 out of 322 scenarios are generated based on 10 basic scenarios by scaling up and down the total demand to study the impact of overall demand on the performance of the proposed model. These scenarios facilitate a comprehensive evaluation of the proposed model's performance under various traffic conditions and demand levels. The through movements at the main intersection consist of two lanes, while the right-turn movements at the main intersection are all single-lane. The approaches from the main intersection to the sub-intersections also feature two lanes. A summary of the input parameters used in the numerical experiments is provided in Table 1. The same set of parameters is used in both the proposed model and the benchmark model to ensure consistency in the comparisons.

**Table 1.** Summary of input model parameters.

| Parameter | Value |
|---|---|
| $\alpha_{l,i}$ | 1 for single-lane road sections; 0.65 for two-lane road sections |
| $s_{l,i}$ | 0.75 vehicle/second |
| $\delta$ | 3 seconds |
| $C_{\min}$ | 40 seconds |
| $C_{\max}$ | 150 seconds |
| $g_{\min,i}$ | 10 seconds for all movements |
| $g_{\max,i}$ | 140 seconds for all movements |
| $\eta_i$ | $Q_i / \sum_i Q_i$ for each movement $i$ |
| $b_e$ | 5 seconds |
| $t_{d,i}$ | 15 seconds for right-turn movements from a sub-intersection to the main intersection; 11 seconds for all movements from the main intersection to sub-intersections; 26 seconds for right-turn movements from one sub-intersection to another sub-intersection |
| $z_{d,i}$ | 9 vehicles for the right-turn bay; 10 vehicles for the bay between sub-intersection and the main-intersection; 30 vehicles for the bay between the main intersection and the entrance of right-turn bay |

*Results and Analysis*

A summary of the results from the conducted simulations to evaluate the performance of the proposed model versus the benchmark model is provided in Table 2. We use two main metrics, including average queue length and average delay. On average, the proposed optimization model reduces delay by 17% and queue length by 31%, with smaller or comparable standard deviation compared to the benchmark model.

**Table 2.** Overall comparison between the benchmark and the proposed model.

| Indicator | Average Delay (s/km) | | Average Queue Length (veh) | |
|---|---|---|---|---|
| | Benchmark model | Proposed model | Benchmark model | Proposed model |
| Average | 61.61 | **51.26** | 173.09 | **117.82** |
| Stdv. | 18.92 | **15.77** | 71.05 | **54.08** |

The simulated and evaluated scenarios are grouped into balanced and imbalanced demand scenarios. A balanced demand represents a scenario in which the sum of the traffic demand of movements 1, 2, 5, and 6 is equal to the sum of the traffic demand of movements 3, 4, 7, and 8. See Fig. 3(a) for an example. Meanwhile, an imbalanced demand represents a scenario in which the sum of the traffic demand of movements 1, 2, 5, and 6 is significantly higher or lower than that of movements 3, 4, 7, and 8. See Fig. 3(b) for an example. The average delay and average queue length for both balanced and imbalanced scenarios are comparatively summarized in Table 3. Results demonstrate that our proposed model outperforms the benchmark model for both balanced and imbalanced scenarios.

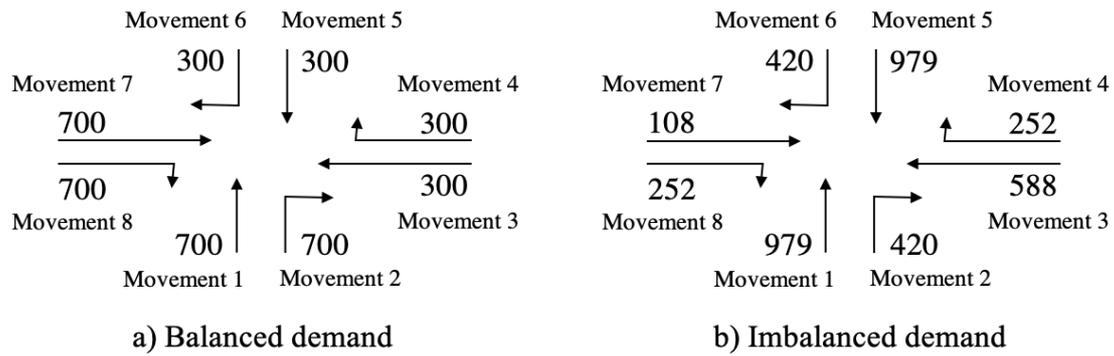

a) Balanced demand          b) Imbalanced demand

**Fig. 3.** Examples of a) balanced and b) imbalanced demand scenarios.

**Table 3.** Results of Balanced Scenarios and Imbalanced Scenarios.

| Scenario | Indicator | Average Delay (s) | | Average Queue Length (veh) | |
|---|---|---|---|---|---|
| | | Benchmark model | Proposed model | Benchmark model | Proposed model |
| Balanced Scenarios | Average | 58.27 | **47.79** | 161.41 | **109.97** |
| | Stdv. | 16.11 | **14.62** | 53.25 | **51.98** |
| Imbalanced Scenarios | Average | 65.87 | **55.80** | 188.03 | **127.88** |
| | Stdv. | 21.87 | **16.43** | 88.93 | **56.03** |

The average delay and average queue length profiles of the 10 base scenarios with different demand levels are illustrated in Fig. 4 to Fig. 6. For each of the 10 base scenarios, a wide range of total demand levels were tested, covering from a very low level of demand to near-maximum demand levels without oversaturating the intersection.

Different levels of total demand in a scenario are generated by scaling up and down the base scenario. Since the traffic demand distributed to the 8 movements is different in each of the base scenarios, slight oversaturation may still occur on different links for different base scenarios when scaling up the demand. As a result, the maximum total demand that keeps the intersection undersaturated differs in each scenario, which can be observed in Fig. 4 to Fig. 6.

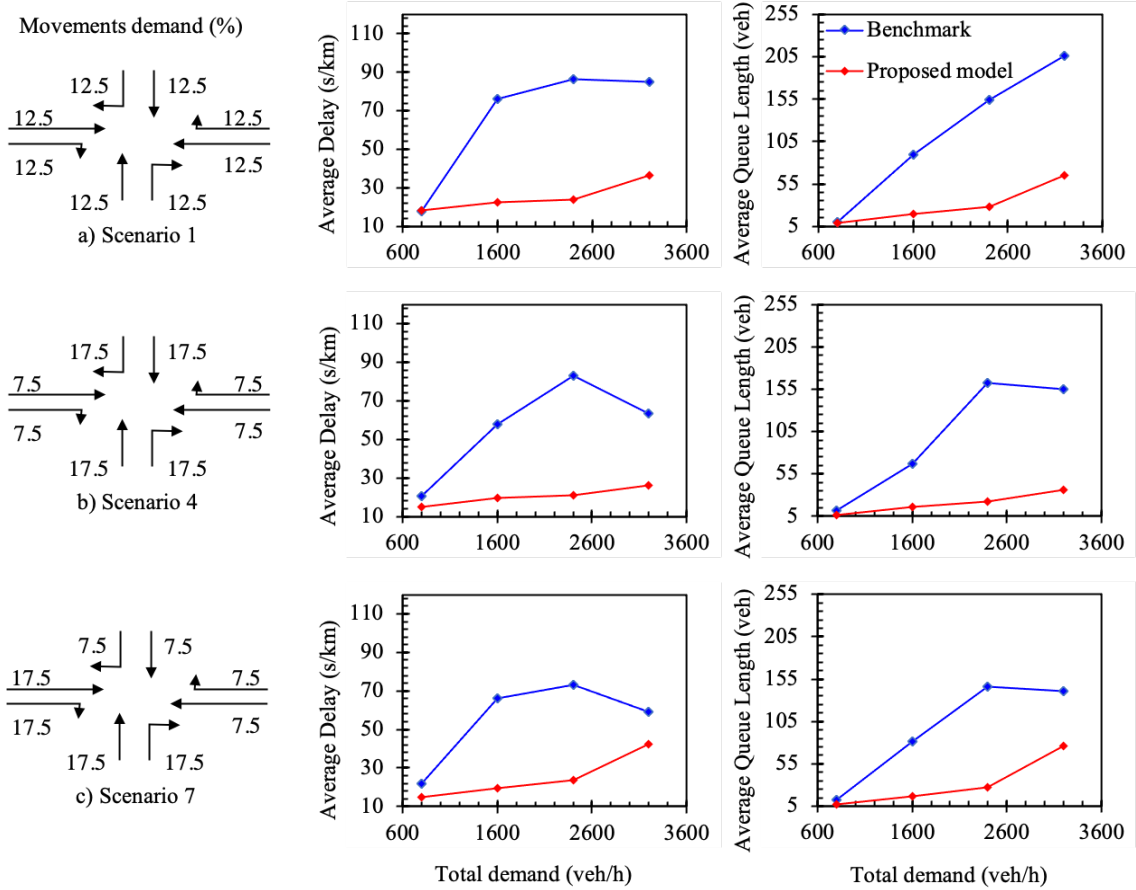

**Fig. 4.** Average delay and average queue length of scenarios with equal demand for through and right-turn movements, a) scenario 1, b) scenario 4 and c) scenario 7.

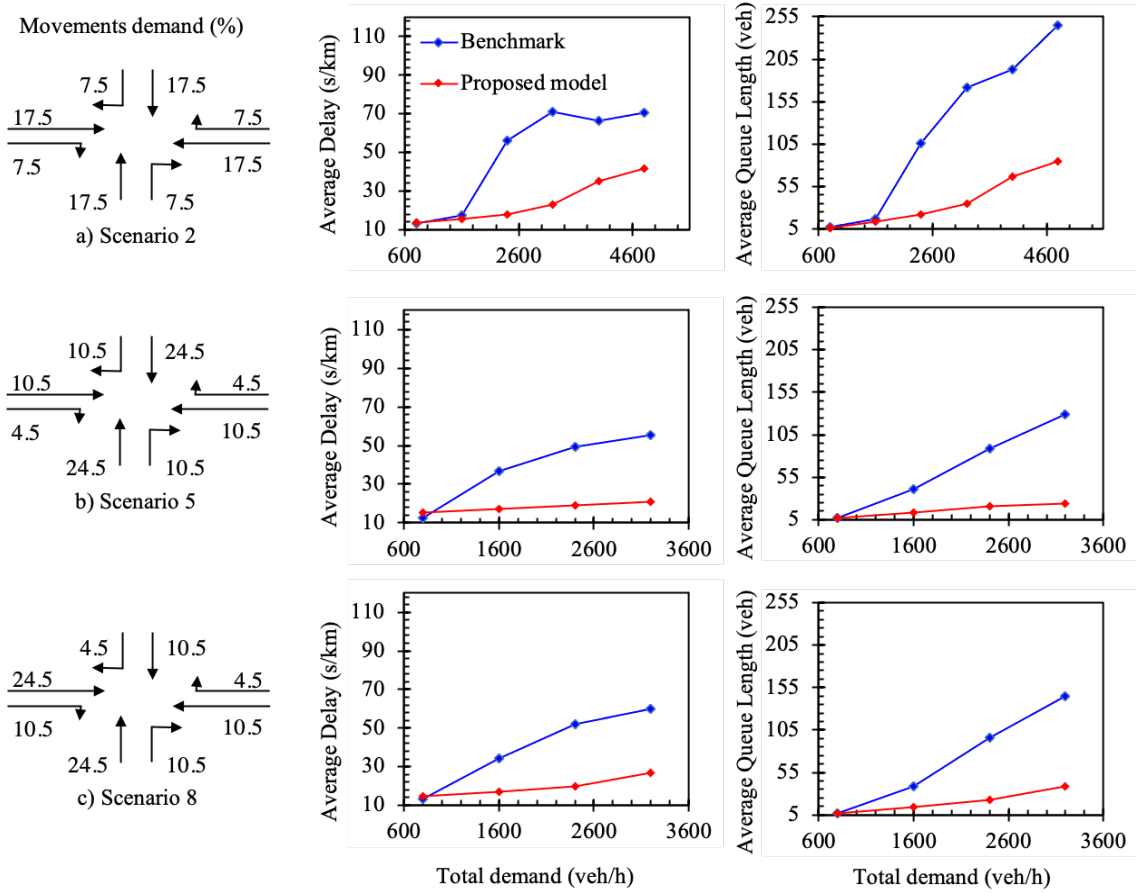

**Fig. 5.** Average delay and average queue length of scenarios with high through demand a) scenario 2, b) scenario 5 and c) scenario 8.

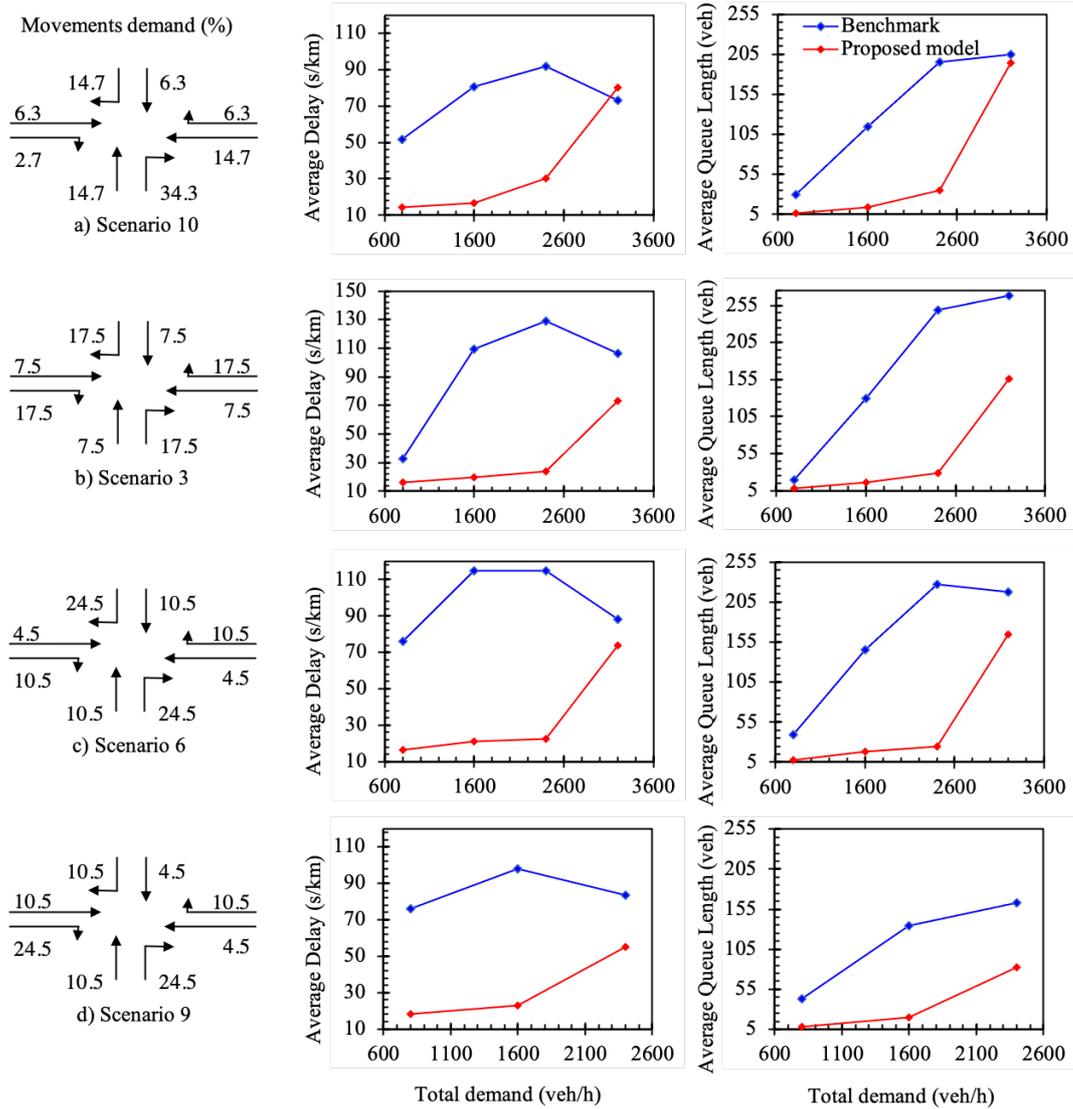

**Fig. 6.** Average delay and average queue length of scenarios with high right-turn demand, a) scenario 10, b) scenario 3, c) scenario 6 and d) scenario 9.

In Scenarios 1, 4, and 7, where the percentage of traffic demand assigned to through movements equals the demand assigned to right-turn movements (see Fig. 4), the proposed model outperforms the benchmark model. In Scenarios 2, 5, and 8, where the percentage of traffic demand assigned to through movements is higher than that assigned to right-turn movements (see Fig. 5), the proposed model performs equally well as the benchmark model in low demand and better in relatively high demand. In Scenarios 3, 6, 9, and 10, where the percentage of demand assigned to through movements is lower than that assigned to right-turn movements (see Fig. 6), the proposed model significantly outperforms the benchmark model. CFIs primarily aim to improve traffic efficiency and safety at intersections with heavy right-turn traffic. Our simulation results suggest that the proposed model creates more efficient signal timing plans than the benchmark model for intersections with heavy right-turn traffic, regardless of the total demand. In some scenarios, the delay of the benchmark decreases as the total demand increases. This is because when demand is high, spillback occurs in some sections, and some vehicles cannot enter the network. The simulator does not include the delay of those vehicles. Such a decrease in demand is not observed in the corresponding results of the proposed model, indicating that compared to the benchmark, the proposed model can avoid spillback when demand is high.

**Table 4.** Summary results for CFI right-turn movements.

| Movement | Indicator | Average Delay (sec/km) | | Average demand (veh) |
|---|---|---|---|---|
| | | Benchmark model | Proposed model | |
| Movement 2 | Average | 166.98 | **163.52** | 484 |
| | Stdv. | 132.57 | **109.12** | |
| Movement 4 | Average | 133.83 | **95.30** | 420 |
| | Stdv. | 108.43 | **24.51** | |
| Movement 6 | Stdv. | 219.81 | **155.10** | 580 |
| | Stdv. | 138.65 | **70.31** | |
| Movement 8 | Average | **66.13** | 78.89 | 332 |
| | Stdv. | 74.69 | **19.21** | |

To further analyze the performance of the proposed model at the movement level, the four right-turn movements across all 322 scenarios are evaluated. See Table 4 for a summary of the analysis outcomes. The proposed model performs significantly better than the benchmark model for movements 4 and 6, while movement 2 exhibits similar performance compared to the benchmark model. The overall statistics in Table 2 suggest that, on average, the proposed model achieves lower levels of delay than the benchmark model. However, for movement 8, which has the lowest level of demand among the four right-turn movements, the proposed model exhibits a larger average delay compared to the benchmark model. These results indicate that the benchmark model allocates unnecessarily long green time to movements with low traffic demand, which affects the efficiency of the entire CFI. Time-dependent average delay profiles for two sample movements are also illustrated in Fig. 7. We compare the performance of the proposed model against the benchmark model for both balanced and imbalanced demand scenarios while keeping the total demand identical for right-turns and through movements. In both cases, the proposed model consistently outperforms the benchmark model.

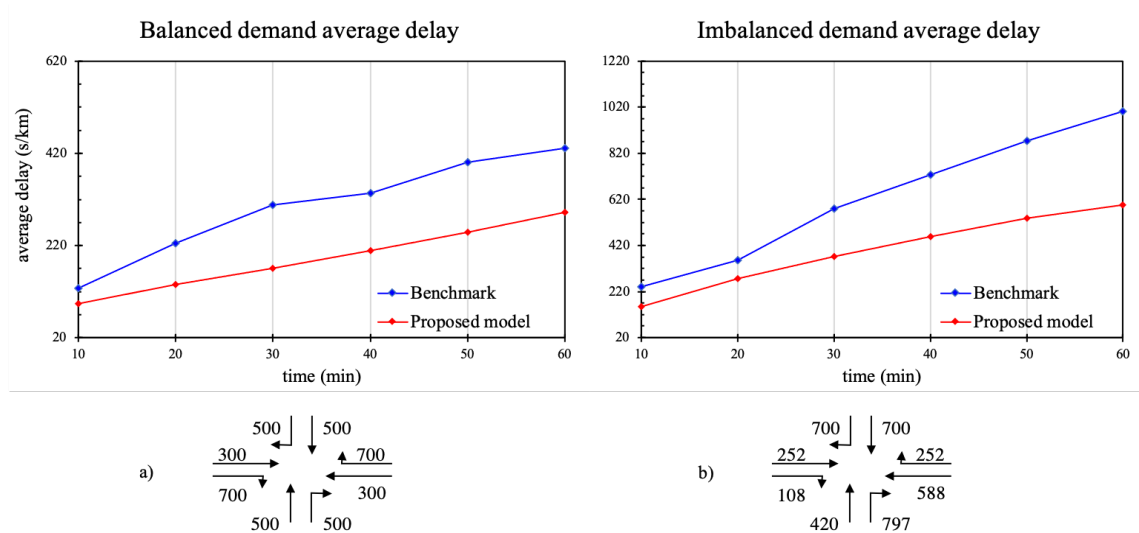

**Fig. 7.** Examples of time-dependant average delay profiles for a a) balanced and b) imbalanced demand scenario.

**CONCLUSION**

In this study, a traffic signal timing optimization model is developed for a full Continuous Flow Intersection (CFI). The performance of the model is tested and compared against a state-of-the-practice benchmark model using over 300 different demand scenarios. The simulation results show that the proposed model improves the traffic efficiency of the intersection compared to the benchmark model for both balanced and imbalanced demand scenarios, regardless of the level of demand. The unique design of CFIs is intended to reduce the impact of heavy right-turn traffic on overall intersection efficiency. The proposed model significantly outperforms the benchmark model in scenarios with high right-turn demand, with smaller average delay by 17% and smaller average queue length by 31%. Future research can explore the inclusion of pedestrians and cyclists in the CFI intersection geometric and traffic signal design and optimization.

**Appendix A: Constraints for $\phi_{l,i}$**

Additional constraints for $\phi_{l,i}$ are listed as below:

$$\phi_{1,5} = \phi_{1,8} \tag{25}$$
$$\phi_{2,7} = \phi_{2,2} \tag{26}$$
$$\phi_{3,1} = \phi_{3,4} \tag{27}$$
$$\phi_{4,3} = \phi_{4,6} \tag{28}$$
$$\phi_{5,1} = \phi_{5,2} = \phi_{5,5} = \phi_{5,6} \tag{29}$$
$$\phi_{5,3} = \phi_{5,4} = \phi_{5,7} = \phi_{5,8} \tag{30}$$
$$\phi_{1,2} + \phi_{1,5} = \phi_{2,4} + \phi_{2,7} = \phi_{3,6} + \phi_{3,1} = \phi_{4,8} + \phi_{4,3} = \phi_{5,1} + \phi_{5,3} = 1 \tag{31}$$

## Appendix B: Constraints for $r_{l,i}$

Additional constraints for $r_{l,i}$ are listed as below:

$$r_{1,5} = r_{1,8} \tag{32}$$
$$r_{2,7} = r_{2,2} \tag{33}$$
$$r_{3,1} = r_{3,4} \tag{34}$$
$$r_{4,3} = r_{4,6} \tag{35}$$
$$r_{5,1} = r_{5,2} = r_{5,5} = r_{5,6} \tag{36}$$
$$r_{5,3} = r_{5,4} = r_{5,7} = r_{5,8} \tag{37}$$
$$r_{1,2} + r_{1,5} = r_{2,4} + r_{2,7} = r_{3,6} + r_{3,1} = r_{4,8} + r_{4,3} = r_{5,1} + r_{5,3} = 1 \tag{38}$$

## Appendix C: Constraints for defining green band cases

A binary variable $y_{d,i,k}$ is used to ensure only one of the six cases are assigned to each movement (see Equation (14)). For case 1, Equation (39), (40), (43) and (47) should be active. For case 2, Equation (42), (43), (44) and (45) should be active. For case 3, Equation (42) and (46) should be active. For case 4, Equation (40) and (44) should be active. And for case 5, Equation (41) and (43) should be active.

$$\theta_a + r_{a,i} + t_{d,i} + n_{a,i}C \geq \theta_b + r_{b,i} + g_{b,i} + n_{b,i}C - C - y_{d,i,k}M$$
$$\forall d \in D_i, (a,b) \in d, i \in I, k = 1 \quad (39)$$

$$\theta_a + r_{a,i} + t_{d,i} + n_{a,i}C \leq \theta_b + r_{b,i} + n_{b,i}C + y_{d,i,k}M$$
$$\forall d \in D_i, (a,b) \in d, i \in I, k = 1,4 \quad (40)$$

$$\theta_a + r_{a,i} + t_{d,i} + n_{a,i}C \geq \theta_b + r_{b,i} + n_{b,i}C - y_{d,i,k}M$$
$$\forall d \in D_i, (a,b) \in d, i \in I, k = 2,5 \quad (41)$$

$$\theta_a + r_{a,i} + t_{d,i} + n_{a,i}C \leq \theta_b + r_{b,i} + g_{b,i} + n_{b,i}C + y_{d,i,k}M$$
$$\forall d \in D_i, (a,b) \in d, i \in I, k = 2,3 \quad (42)$$

$$\theta_a + r_{a,i} + g_{a,i} + t_{d,i} + n_{a,i}C \leq \theta_b + r_{b,i} + g_{b,i} + n_{b,i}C + y_{d,i,k}M$$
$$\forall d \in D_i, (a,b) \in d, i \in I, k = 1,5 \quad (43)$$

$$\theta_a + r_{a,i} + g_{a,i} + t_{d,i} + n_{a,i}C \geq \theta_b + r_{b,i} + g_{b,i} + n_{b,i}C - y_{d,i,k}M$$
$$\forall d \in D_i, (a,b) \in d, i \in I, k = 2,4 \quad (44)$$

$$\theta_a + r_{a,i} + g_{a,i} + t_{d,i} + n_{a,i}C \leq \theta_b + r_{b,i} + n_{b,i}C + C + y_{d,i,k}M$$
$$\forall d \in D_i, (a,b) \in d, i \in I, k = 2 \quad (45)$$

$$\theta_a + r_{a,i} + g_{a,i} + t_{d,i} + n_{a,i}C \geq \theta_b + r_{b,i} + n_{b,i}C + C - y_{d,i,k}M$$
$$\forall d \in D_i, (a,b) \in d, i \in I, k = 3 \quad (46)$$

$$\theta_a + r_{a,i} + g_{a,i} + t_{d,i} + n_{a,i}C \geq \theta_b + r_{b,i} + n_{b,i}C - y_{d,i,k}M$$
$$\forall d \in D_i, (a,b) \in d, i \in I, k = 1 \quad (47)$$

**Data Availability Statement**

Some data and models that support the findings of this study are available from the corresponding author upon reasonable request. The developed code is, however, proprietary to the funding organizations.

**Notation List**

The following symbols are used in this paper:

| Sets | |
|---|---|
| $L = \{1, 2, 3, 4, 5\}$ | Set of intersections (see Figure1, Intersection 1 represents the south sub-intersection, intersection 2 represents the east sub-intersection, intersection 3 represents the north sub-intersection, intersection 4 represents the west sub-intersection, and intersection 5 represents the main intersection) |
| $L_i$ | Set of intersections passed by movement $i$ |
| $I = \{1, \ldots, 8\}$ | Set of traffic movements at the main intersection (index of each movement is shown in Fig. 1) |
| $K$ | Set of indices for uncoordinated time periods |
| $D_i$ | Set of pairs of intersections that passed by route $i$ in sequence |
| **Parameters** | |
| $Q_i$ | Average traffic demand for movement $i$ (per second) |
| $\alpha_{l,i}$ | Lane use factor at intersection $l$ for movement $i$ |
| $s_{l,i}$ | Saturation flow rate per lane at intersection $l$ for movement $i$ (per second) |
| $\delta$ | Lost time (seconds) |
| $C_{\min}$ | Minimum cycle length (seconds) |
| $C_{\max}$ | Maximum cycle length (seconds) |
| $g_{\min,i}$ | Minimum green time for movement $i$ (seconds) |
| $g_{\max,i}$ | Maximum green time for movement $i$ (seconds) |
| $\eta_i$ | Weighting factor of movement $i$ |
| $b_e$ | Minimum effective green bandwidth (seconds) |
| $t_{d,i}$ | Travel times used in movement $i$ between intersection pair $d$ (seconds) |
| $z_{d,i}$ | Capacity of the bay between intersection pair $d$ for movement $i$, $d = 0$ when representing the maximum acceptable distance between stop line and end of queue at the first intersection passed by movement $i$ (vehicles) |
| $M$ | A very large positive value |
| **Variables** | |
| $\mu_l$ | A multiplier for the traffic volume at intersection $l$, indicating the ratio between actual capacity and the volume, $\in [0,1]$ |
| $\xi$ | The reciprocal of the common cycle length, $\in [1/C_{\min}, 1/C_{\max}]$ |
| $C$ | Common cycle length in seconds, $\in \mathbb{Z}$ |
| $\phi_{l,i}$ | Green time duration at the intersection $l$ for movement $i$ expressed as the proportion of the cycle time, $\in (0,1)$ |
| $g_{l,i}$ | Green time at intersection $l$ for movement $i$, $\in \mathbb{Z}$ |
| $\tau_i$ | Maximum distance between stop line and rear of queue tolerated for each direction in number of vehicles, $\in \mathbb{Z}$ |
| $b_i$ | Bandwidth for movement $i$, $\in \mathbb{Z}$ |
| $\theta_l$ | Offset of intersection $l$, $\in \mathbb{Z}$ |
| $r_{l,i}$ | Time between the start of a cycle and the start of the green phase used by movement $i$ at intersection $l$ (seconds), $\in \mathbb{Z}$ |
| $w_{l,i}$ | Time between the start of green and the start of band for movement $i$ at intersection $l$, $\in \mathbb{Z}$ |
| $n_{li}$ | An integer indicator of cycle length, indicating the start of bandwidth $\in \{0,1\}$ |
| $y_{d,i,k}$ | Variable for selecting green-band case $k$ for movement $i$ between intersection pair $d$, defined in Equation (18), $\in \{0,1\}$ |
| $P_{d,i,k}$ | Uncoordinated time period of movement $i$ between intersection pair $d$ for case $k$, $\in \mathbb{Z}$ |

**Author contributions**
The authors confirm contributions to the paper are as follows: study conception and design: Y. Hu, M. Saberi, R. Mohajerpour and D. Rey; model development: Y. Hu, D. Rey, R. Mohajerpour, M. Saberi; optimisation model coding: Y. Hu; analysis and interpretation of results: Y. Hu, M. Saberi and D. Rey; draft manuscript preparation: Y. Hu, D. Rey, R. Mohajerpour and M. Saberi. All authors reviewed the results and approved the final version of the manuscript.

**Acknowledgments**
We are thankful to iMOVE CRC and Synergistic Traffic Consulting for sponsoring this research project. We would also like to acknowledge Johnny Leung and Inv Valiant Yuk Yuen Leung from Synergistic Traffic Consulting for their support. We are also grateful to Saeed Rahmani for his contributions to the project in its early stages.

**Funding**
This research has been funded by iMOVE CRC and Synergistic Traffic Consulting.